\documentclass[12pt,a4paper]{article}
\usepackage{graphicx}
\usepackage{newtxtext}
\usepackage{newtxmath}
\usepackage{hyperref}

\newtheorem{theorem}{Theorem}
\newtheorem{lemma}{Lemma}

\newcommand{\RomanNumeralCaps}[1]
\title{Existence and Smoothness of the Navier-Stokes equation \\
  using the \\
  Boundary Integral Method}
\author{Edmund Chadwick, University of Salford}
\date{e.a.chadwick@salford.ac.uk}

\begin{document}

\maketitle

%\shorttitle{Existence, Smoothness, Uniqueness and Validation of the
%  Navier-Stokes equation} 

%\author{Edmund Chadwick}
%  \corresp{\email{e.a.chadwick@salford.ac.uk}}}
      
%\author{Edmund Chadwick\aff{1}
%  \corresp{\email{e.a.chadwick@salford.ac.uk}}}

%\affiliation{\aff{1}Mathematics, School of Computing, Science and
%  Engineering, University of Salford, Salford M5 4WT,
%  UK\email{$^*$Corresponding author: Edmund Chadwick }} 

%\author[1]{\fnm{Edmund} \sur{Chadwick}}\email{e.a.chadwick@salford.ac.uk}
%\affil[1]{\orgdiv{Mathematics}, \orgname{University of Salford}, 
%  \orgaddress{\street{} \city{Salford} \postcode{M54WT},
%    \state{} \country{U.K}}} 

\begin{abstract}
  %\abstract{
  Consider an exterior space-time domain where the incompressible
  Navier-Stokes equation and continuity equation hold with no bodies
  or force fields present, 
  and smooth velocity at initial time. 
  A smooth solution with a stokeslet far-field decay for all
  subsequent time is sought and found,
  demonstrating existence and smoothness.
  A space-time boundary integral velocity representation is given by
  an integral distribution of fundamental solutions of the
  Navier-Stokes equation called nslets.
  These nslets approach eulerlets close to their origin which have a
  singularity line in the fluid that moves with the fluid to ensure
  that the velocity direction is defined.
  The boundary enclosing the fluid point is chosen to move with the
  fluid also and so in this reference frame the Lagrangian material
  derivative and Eulerian partial derivative become the same in the
  limit.  
  Consequently, the contributions to the flux from the quadratic terms 
  originating from the non-linearity vanish thereby enabling the
  boundary integral method standard theory of Oseen and Ladyzhenskaya
  to be used for this non-linear problem. 
  It is then shown that the resulting representation exists and is
  smooth. 
  Zero initial velocity gives the null solution.
  The non-linear interaction between the flow field and the
  fundamental solution alignment to it describes a dynamical system
  of two interacting linear systems incorporating chaos,
  and an example demonstrating reduction to the blinking vortex is
  given.   
\end{abstract}
%}

%\begin{keywords}
%Navier-Stokes equations, Green's Boundary Integral Method
%\end{keywords}

%\keywords{Navier-Stokes equations, Green's Boundary Integral Method,
%  tensor analysis}

%\maketitle

\section{Introduction}

\subsection{The key idea}
The key idea is that Lagrangian and Eulerian descriptions coincide for
a space-time boundary integral in the limit approaching and enclosing
a fluid point and moving with the fluid,
and consequently the non-linearity disappears. 
This is specifically the boundary integral contribution in the
boundary integral method used for fluids by Oseen \cite{Oseen:1927}
and Ladyzhenskaya \cite{Ladyzhenskaya:1969}.
This results in a boundary integral representation for the velocity
from which existence and smoothness are determined.

\subsection{Millennium problem}

A general solution to the Navier-Stokes equation is sought which is
one of the millennium problem challenges \cite{Fefferman:2000}.
This challenge seeks to prove (or otherwise) that there exists a
smooth solution satisfying the Navier-Stokes equation for the particular
problem of a smooth divergence free velocity at initial time and an
incompressible fluid filling all of space.  
Existence and smoothness have been proven for the two-dimensional case,
and is demonstrated by Ladyzhenskaya \cite{Ladyzhenskaya:1969} for the weak
solution by using generalised functions.
However, this approach has not successfully been extended to the
three-dimensional unsteady problem. 
This is because of the question of blow-up of the solution due to a
finite-time singularity;  
At small scale the inertial (transport) term dominates and the
Navier-Stokes equation becomes supercritical and has been demonstrated
to blow up when considering averaged quantities \cite{Tao:2016}.
For the three-dimensional steady Stokes problem,
Ladyzhenskaya also gives a strong solution by a different approach,
using the theory of hydrodynamic potentials which uses Green's
boundary integral method.
Oseen \cite{Oseen:1927} similarly gives boundary integral solutions
for the three dimensional steady Oseen flow,
and also unsteady Stokes and Oseen flow using stokeslet and oseenlet
unsteady force impulse fundamental solutions. 
Ramm \cite{Ramm:2019} gives a velocity representation for the
unsteady Navier-Stokes problem by a domain distribution of unsteady
stokeslets, 
with the non-linear velocity term within the domain integral itself.
By using Fourier transforms,
he then is able to give bounds for the velocity and velocity
derivative as well as demonstrate the solution is unique.
However, it is unclear how this approach deals with inertial blow-up
singularity and there has been a challenge to the accuracy of the
analysis in the review \cite{LemarieRieusset:2019} that claims one of
the integrals used in the analysis is logarithmic unbounded. 

\subsection{Boundary integral nslet representations} 

Chadwick \cite{Chadwick:2018a} \cite{Chadwick:2019a}
\cite{Chadwick:2023} also considers fundamental solutions which are
nslets for the Navier-Stokes equations that have the expected singular
behaviour, rather than stokeslets that don't, 
and develops a boundary integral representation from the formulations
used in fluid dynamics by Oseen and Ladyzhenskaya. 
Chadwick first develops this approach giving a boundary integral
representation for Euler flow \cite{Chadwick:2019} and then
subsequently Navier-Stokes flow \cite{Chadwick:2018a}
\cite{Chadwick:2019a} by introducing their fundamental solutions.
From this a formulation is given \cite{Chadwick:2023},
such that nslets reduce to supercritical inertia dominant eulerlets at
small scale near to their origin.
Although the eulerlets are singular in this domain,
their integral flux is finite and so can be evaluated in the boundary
integral representation giving finite fluid velocity (and pressure).
So although the eulerlet contains a singularity,
the integral flux does not enabling the boundary integral
representations \cite{Chadwick:2019a} \cite{Chadwick:2023}. 
This formulation has been successfully tested against benchmark
problems discussed next. 

\subsection{Benchmark tests for the boundary integral representation}

For low Reynolds number two-dimensional steady flow past a
circular cylinder, the leading order approximation to the nslet is the
two-dimensional steady oseenlet.
So this provides a benchmark test, since this approximation is
expected to give good results according to the nslet representation of
the new theory but poor results according to standard theory because
the Oseen approximation violates the boundary condition.
This benchmark test proved a successful validation of the new theory
where good agreement is given for velocity,
pressure and drag coefficient that wouldn't have been expected from
standard theory
\cite{Bwebum:2019} \cite{Bwebum:2020} \cite{Bwebum:2020a}
\cite{Bwebum:2025}
.  
Even for Reynolds number greater than one the resulting eddy
description is modelled closely \cite{Bwebum:2025}.

High Reynolds number non-linear steady flow past a flat plate has been
considered. 
For a semi-infinite flate plate the boundary layer is satisfied by the
Blasius non-linear equation.
The new formulation gives a representation that accurately models the
Blasius solution outperforming all other approaches
%such as the
%shooting method, and (Chebyshev) polynomial representations
\cite{AdamuChadwick:2023}. 
For a finite flat plate the solution is given by triple deck matched
asymptotic theory.
The new formulation gives the Goldstein near-wake formulation and all
the expected matched asymptotic triple deck relations
\cite{Chadwick:2025}. 
This results in an inner deck boundary element formulation
%giving the expected velocity and pressure variation
\cite{Chadwick:2026}.
So the new theory accurately models these two flat plate non-linear
benchmark problems.

\subsection{New developments in the present paper}

In the present paper, formulations are given with greater detail
particularly for the boundary integral calculations,
and as well correcting for the singularity velocity direction of the
small scale eulerlet limit which is set to move along a fluid
streakline ensuring continuity in the velocity vector direction.
This is a crucial inclusion, as this admits chaotic and turbulent
representations, and consequently all the evaluations have been
recalculated in the present paper in the reference frame moving
with the velocity at the fluid point.  
Consequently, the space-time boundary enclosing the fluid point is
aligned to this moving reference frame.
The boundary integral methodology of Oseen \cite{Oseen:1927} is
followed that uses the Green's integral representation theory of
hydrodynamic potentials for the linear problems of Oseen and Stokes
flow but extended by considering a four-dimensional space-time domain
with boundaries on space-time hypersurfaces.
As well, a divergence function representing non-linear contributions
and directed along the fundamental solution radius is considered in
order to then be able to apply the divergence theorem to get a
boundary integral formulation.  
This has similarities to the dual reciprocity method especially in
that the non-linear term relies on evaluations within the domain
itself. 
The resulting integral evaluation from the non-linear divergence is
then shown to give no contribution to the Green's boundary integral
formulation for the velocity over the far-field boundary and give a
contribution that vanishes for the near-field boundary
enclosing the fluid point in the limit as this boundary size tends to
zero.
%So the non-linear term drops out as expected for space-time boundaries
%that move with the fluid so retaining momentum flux,
%since the origin of the non-linear term in the Eulerian description of
%the Navier-Stokes equation is from the momentum flux leaving a control
%volume. 
The velocity is then represented by a linear distribution of nslet
fundamental solutions, and existence and smoothness follow from this
description.

The centres of the fundamental solutions themselves are directed along
the fluid velocity direction and so depend upon the velocity at that
point.
This is a dynamical system that gives rise to the possibility of chaos
and turbulence. 
For example,
for high Reynolds number limit Euler flow for a set of blinking
vortices,
the eulerlet singularity direction using the new theory is given by
the velocity flow field at that point. 
Assuming leading order steady two-dimensional flow,
the new formulation then trivially reduces to an equivalent
representation as that for standard blinking vortex chaotic mixing,
see for example \cite{DaitcheTel:2009}.

Therefore, the new representation has been shown to successfully model
low Reynolds number flow, high Reynolds number non-linear flow and
chaotic high Reynolds number limit Euler flow. 
The final benchmark currently being investigated is to model
turbulence.

\subsection{Previous work by the author on the foundations for the new
  formulation} 

This boundary integral approach builds upon previous work on providing
boundary integral formulations in incompressible fluid dynamics by the
author which has incrementally developed Oseen's mathematical analysis
for boundary integral representations and is outlined next.
The work started by considering the far-field Oseen representation
for large Reynolds number \cite{Chadwick:1992},
which resulted in a far-field velocity representation by an integral
distribution of oseenlets \cite{Chadwick:1998}
\cite{FishwickChadwick:2006}. 
From this, a slender body theory in Oseen flow was developed
\cite{Chadwick:2002}, and verified experimentally in low-speed wind
tunnels \cite{Chadwick:2010}.
In the large Reynolds number limit,
this then gives insight into the form of the near-field potential
flow, and a slender wing theory in potential flow was developed
\cite{Chadwick:2005} as well as an understanding of the underlying
vortex flow structure \cite{Chadwick:2006} \cite{Chadwick:2007}.
The work was then recast within a matched asymptotic framework
but using a novel approach with the boundary integrals of the
representations being matched.
This was achieved at low Reynolds number for an outer Oseen flow
and an inner Stokes flow \cite{Chadwick:2013} and also later in
\cite{Bwebum:2019} \cite{Bwebum:2020} \cite{Bwebum:2020a}. 
The method was also applied at high Reynolds number to matching 
an outer Oseen flow to an inner Euler flow
\cite{Kapoulas:2013}
\cite{Chadwick:2015} \cite{Chadwick:2015a}
\cite{Chadwick:2018} \cite{Chadwick:2019}.
This leads to a velocity representation in Euler flow by an integral
distribution of Euler flow fundamental solutions called eulerlets.
As well, it is shown that the representation reduces to all the
standard classic steady Euler flow models (potential aerofoil theory, 
vortex-lattice theory, thin-wing theory, slender-wing theory,
slender-body theory) under the appropriate approximations
\cite{Chadwick:2019}.  
Furthermore, this theory was able to model bluff body flows by an
Euler slip wake expected from the theory \cite{Chadwick:2019}.

\section{Problem statement}
This paper seeks to prove the following theorem.

\begin{theorem}
  Consider the motion of a fluid in $\mathbb{R}^3$ described by the
  incompressible Navier-Stokes equation and continuity equation,
  such that body forces are taken to be zero, 
  for the unknown velocity field $\in \mathbb{R}^3$ and pressure $\in 
  \mathbb{R}$ defined for time $t \ge 0$,
  with initial condition given for the velocity as a $C^{\infty}$
  divergence free vector field that decays like a stokeslet in the
  far-field.
  Then, the pressure and velocity exist and are smooth at all
  subsequent time with bounded energy.
\end{theorem}

This theorem requires the following two lemmas.

\begin{lemma}
The problem can be restated as a boundary value problem in a space-time
domain $\Sigma _+$ enclosed by a hypersurface $\partial \Sigma$
consisting of a spherinder (spherical cylinder) surface
$\partial \Sigma_X$ given by:
$x_1^2+x_2^2+x_3^2=X^2$ and $ 0 \le t \le X$;
and two volumes $x_1^2+x_2^2+x_3^2 \le X^2$ one at initial time $t=0$,
and the other at time $t=X$;
where ${\bf x}$ is co-ordinate position and $t$ is time,
and $X \rightarrow \infty$, see figure \ref{IVPdomain}. 

\begin{figure}[h]
\begin{center}
\includegraphics[width=9.5cm]{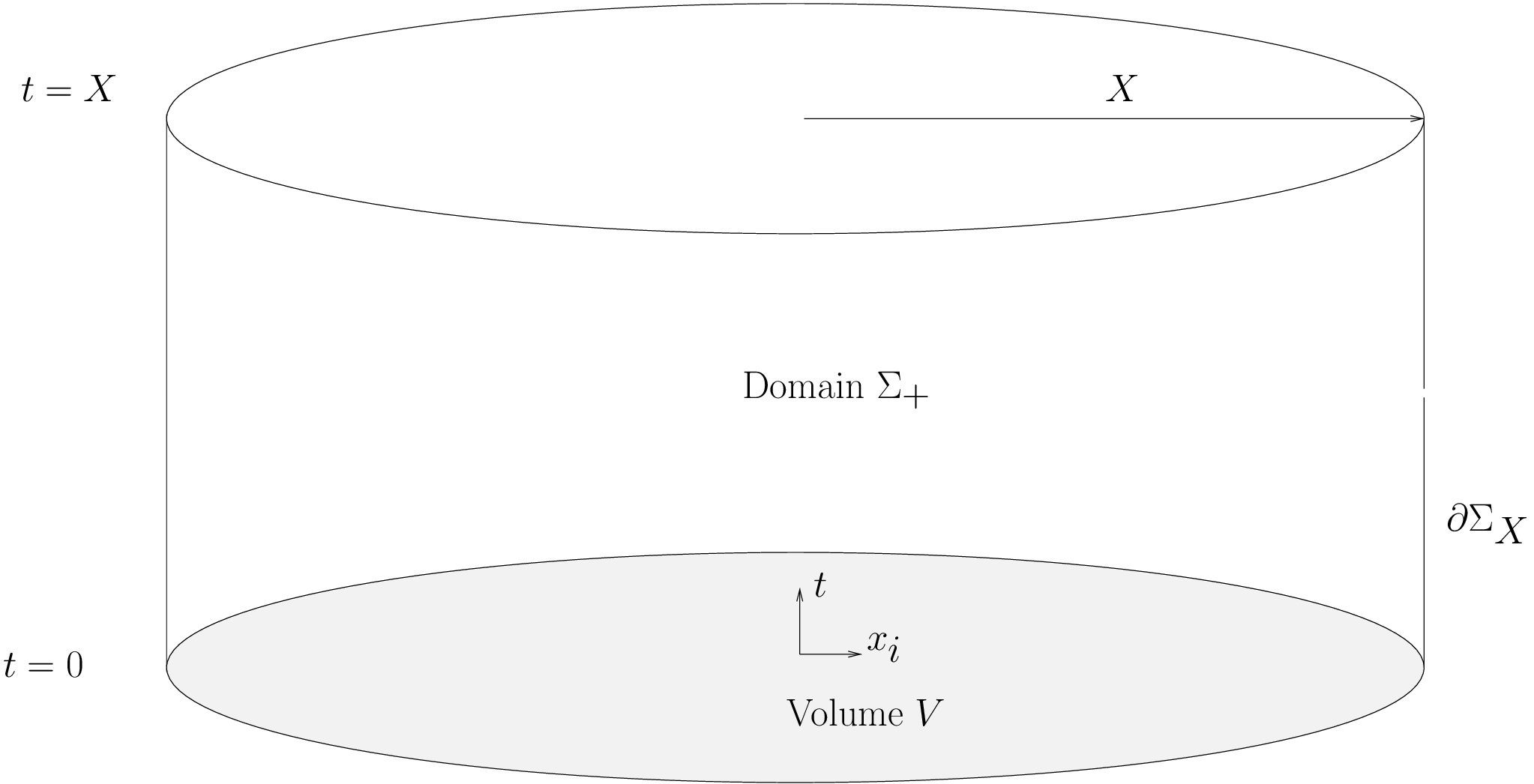}
\end{center}
\caption{The space-time domain $\Sigma ^+$}   
\label{IVPdomain}
\end{figure}

The boundary conditions are:
The fluid velocity stipulated over the hypersurface volume $V$ at
initial time,  
and over the far-field spherinder hypersurface $\partial \Sigma _X$
the flow field behaves like a stokeslet.

The spherinder hypersurfaces have defined normals (except in the limit
approaching the edges), are smooth surfaces piece-wise connected,
and so amenable to Green's boundary integral analysis.
Green's analysis is performed in the space-time domain enclosed by the
initial time and far-field space-time boundary integrals,
and the vanishingly small space-time boundary aligned to the flow
direction and enclosing the space-time point.  
A proof by ansatz is considered.
So, in this domain it is assumed that the velocity and pressure (and
therefore the nslet velocity and pressure also) are $C^{\infty}$
smooth. 
This enables a solution for the velocity to be found by a Green's
integral representation.
From this solution it is then shown that a velocity and pressure which
are $C^{\infty}$ smooth exist, completing the proof by ansatz.
%It is then shown that there can be only one solution, uniqueness,
%in the $C^{\infty}$ smooth class.
In particular, the Green's integral representation for the velocity
$u_k({\bf x},t)$ is given by a boundary integral distribution of nslet
fundamental solutions with strength given by (minus) the initial
velocity, distributed over the all of space volume $V$ at the initial
time $t=0$ so

\begin{equation}
u_k({\bf x},t) = -\int _V u_i({\bf x'},0) u_{ki} ({\bf x}-{\bf x'},t)
dV',
\end{equation}
where $u_{ki}$ is the $i$th velocity component of the $k$th nslet and
is chosen to be aligned to the fluid direction at its origin,
and $dV'$ is an element of the volume $V$.
\end{lemma}

%\begin{proof}
{\bf Proof}
Given in Sections 3-7.
%\end{proof}
\
\begin{lemma}
  From this, the velocity is then given by
\begin{equation}
  u_k({\bf x^L},T) = u_k({\bf x^L},0)(1+{\rm O}(\nu \sqrt{T})),
\end{equation}
where ${\bf x^L}$ is the Lagrangian co-ordinate position moving with 
the fluid having viscosity $\nu$, for some small time $\nu \sqrt{T}$.
Once this is established, it follows that the solution exists and is
smooth for all time $t$.
\end{lemma}

%\begin{proof}
{\bf Proof}
Given in Section 8.

\section{Governing Equations}
Consider the incompressible Navier-Stokes equation and continuity
equation given by
\begin{equation}
  \rho f_i= \rho u_{i,0'}+ \rho u_ju_{i,j'}+p^{\dagger}_{,i'}-\mu
  u_{i,j'j'}=0 \;, \; \; u_{i,i'}=0 , 
\label{GE:NavStokesDimensional}
\end{equation}
where $f_i({\bf x'})$ is the body force taken to be zero,
$u_i({\bf x'})$ is the velocity,
$p^{\dagger} ({\bf x'})$ is the pressure, 
$\rho $ is the fluid density and $\mu $ is the dynamic viscosity.
The capital index refers to space-time co-ordinates $x'_I$, $0
\le I \le 3$ with $x'_0$ being the time variable,
and a lower index refers to spatial Cartesian co-ordinates $x'_i$, $1
\le i \le 3$.
The analysis in this paper uses lengthy complex terms including
differentials that intricately mix the physical variables and
fundamental solutions.   
In order to make these expressions sufficiently concise and readable, 
the Einstein summation index convention is used,
and so $a_I b_I = a_0b_0+a_1b_1+a_2b_2+a_3b_3$.
We also borrow the General Relativity derivative representation by a
comma. 
Furthermore, here the prime over the derivative index represents a
derivative with respect to the primed Cartesian co-ordinate $x'_I$
such that 
$\frac{\partial g({\bf x'})}{\partial x'_I} =  g_{,I'}$,
for some function $g$.
The primed co-ordinate will refer to an integrated value over the
boundary integral given later in the Green's Boundary Integral Method.
To reduce the number of parameters in the problem ($\rho$ and $\mu$)
from two to one, 
it is usual to express the quantities in dimensionless form,
giving a single parameter which is the Reynolds number. 
However, for this problem there is no body and therefore no
characteristic length, and so this is not possible. 
To overcome this, Fefferman \cite{Fefferman:2000} considers dividing
by $\rho$ to get
\begin{equation}
  f_i= u_{i,0'}+u_ju_{i,j'}+p_{,i'}-\nu u_{i,j'j'}=0 \;, \; \; u_{i,i'}=0 ,
\label{GE:NavStokes}
\end{equation}
where $p({\bf x'})$ is the pressure for a density of one
(or pressure divided by the density constant),
and now $\nu =\mu /\rho$ which is the kinematic viscosity is the only
parameter in the equation. 
(Comparing to the dimensionless Navier-Stokes equation, then $\nu$ is
in the position of the reciprocal of the Reynolds number with
Bernoulli pressure (ignoring the half) used as the characteristic
pressure, 
which provides an easy way to transform into dimensionless form for
other problems where there is a characteristic body length.)

This way of representing the Navier-Stokes equation is used widely by
those working on the Navier-Stokes millennium problem,
and so it shall also be followed here. 
For ease of use, like Fefferman and others, we call $p$ the pressure
from now on, but to obtain the standard boundary integral formulation
the actual pressure is easily inserted back into the expression from
$p^{\dagger}=\rho p$. 

The millennium problem is for the time dependent problem,
and so to solve this problem we must consider the space-time domain
and the three fundamental solutions to equation (\ref{GE:NavStokes})
$1 \le k \le 3$ which we call nslets, 
that are generated by applying unit impulse forces in the $k$
direction at the point $x_I = x'_I$,
and such that the solutions decay and vanish in the far-field.
Using the Dirac delta function notation, this is described by
\begin{equation}
  f_{ki}= u_{ki,0}+u_{{\bf k}j} u_{{\bf k}i,j}+p_{k,i}-\nu
  u_{ki,jj}= - \delta (x_I-x_I') \delta _{ki}  \;, \; \;\; 
  u_{ki,i}=0 ,
\label{GE:NS}
\end{equation}
where the first index $k$, denotes the $k$th solution.
The comma represents a derivative with respect to the 
Cartesian co-ordinate $x_I$ such that
$\frac{\partial g({\bf x-x'})}{\partial x_I} =  g_{,I}$.

\begin{figure}[h]
\begin{center}
\includegraphics[width=8.cm]{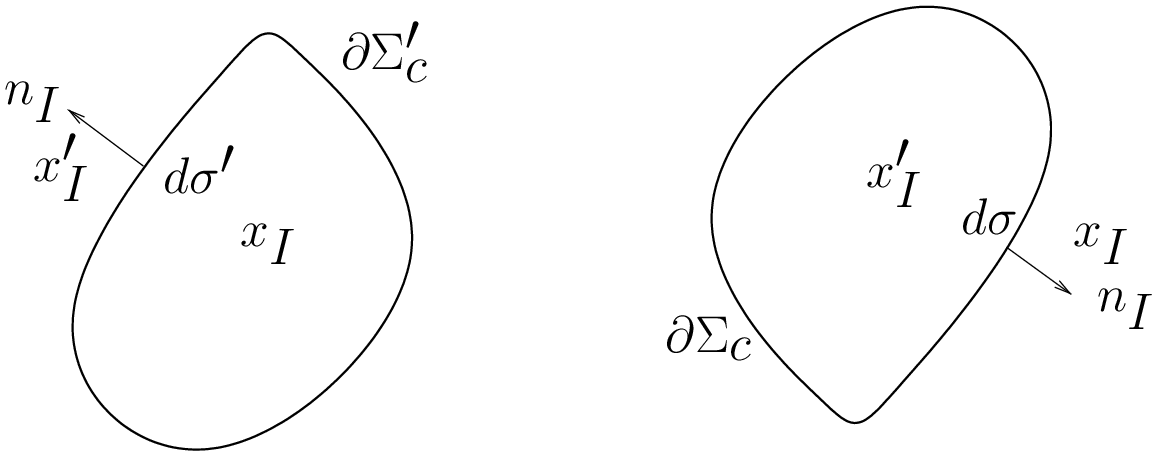}
\end{center}
\caption{The normal direction between integrations for
  the primed and unprimed co-ordinate}
\label{GE:normal}
\end{figure}

A bold index does not use the Einstein summation index
convention, and so $a_{\bf k} b_{\bf k} = a_1b_1$ for $k=1$,
$a_{\bf k} b_{\bf k} = a_2b_2$ for $k=2$ and
$a_{\bf k} b_{\bf k} = a_3b_3$ for $k=3$.
The term $\delta _{ki}$ is Kronecker delta, such that
$\delta _{ki} = 1$ for $i=k$, and 
$\delta _{ki} = 0$ for $i \ne k$.
The Dirac delta function $\delta (x_I-x_I')$ is given in
space-time, 
and so $\delta (x_I-x_I') = \delta (x_0-x_0') \delta (x_1-x_1') \delta
(x_2-x_2') \delta (x_3-x_3')$,
where the Dirac delta function $\delta (x_0-x_0')$ (for example) is
defined as being zero everywhere except at $x_0=x_0'$ such that $\int
_a^b \delta (x_0-x_0') dx_0=1$ if $a  < x_0'$ and $b> x_0'$, and is zero
otherwise.
Hence, the Dirac delta function is a notation and the description
(\ref{GE:NS}) should be read as the equation
\begin{equation}
  \int _{\partial \Sigma _c} f_{kiJ}n_J d\sigma = -\delta _{ki}
\label{GE:fkij-int}
\end{equation}
where:
$\partial \Sigma _c$ is a closed space-time hypersurface boundary
enclosing the point $x'_I$ having an outward pointing normal $n_J$,
and $d\sigma$ is an element of the boundary;
and the function $f_{kiJ}$ is given by
\begin{equation}
  f_{kiJ}=u_{ki}\delta _{0J} + u_{{\bf k}j}u_{{\bf k}i}+p_k\delta
  _{ij}-\nu u_{ki,j},
\end{equation}
where the vector product with mixed space-time and spatial indices is
defined as $a_Jb_j =a_1b_1+a_2b_2+a_3b_3$.
We note that (\ref{GE:fkij-int}) can also be rewritten in terms of an
integral of the primed co-ordinate.
In this case, the unit normal at the point $x_I$ for the surface
integral over elements $d\sigma$ enclosing the primed co-coordinate $x'_I$
is minus the unit normal at the primed co-ordinate point $x'_I$ for
the surface integral over elements $d\sigma'$ enclosing the co-coordinate
$x_I$, see figure \ref{GE:normal}.

\noindent So
\begin{equation}
  n_Jd\sigma ' = -n_Jd\sigma ,
  \label{primetodomain}
\end{equation}
and 
\begin{equation}
  \int _{\partial \Sigma '_c} f_{kiJ}n_J d\sigma ' = \delta _{ki},
\label{GE:fkij-int-primed}
\end{equation}
where $\partial \Sigma '_c$ is the reflection about all four
space-time axes of the boundary $\partial \Sigma _c$, see figure
\ref{GE:normal}. 
It is noted that the solution to (\ref{GE:NS}) is non-unique,
and the particular solution required depends upon the particular flow
conditions at the position in the fluid for the particular problem
under consideration.
So, this gives a different type of boundary integral equation from the
standard whereby now the right hand side kernal (fundamental solution)
is also an unknown as well as the left hand side physical variable
(fluid velocity).
This interaction then under certain circumstances generates chaotic
and turbulent flows. 

\section{NSlet asymptotic approximations}
Existence and smoothness are established from the velocity
representation by nslets.
This representation is obtained from the theory of hydrodynamic
potentials, see for example \cite{Ladyzhenskaya:1969}.
This requires the evaluation of boundary integrals near to and
far from a fluid point, and so the near-field and far-field asymptotic
approximations of the nslet.
So in this section, consider approximations to the nslet near to and
far from its origin. 
It shall be shown that near to its origin, the nslet approximates to
the eulerlet, and far from its origin, the nslet approximates to the
stokeslet. 

\subsection{Near to the origin: the eulerlet}\label{section4a}

Consider the orders of terms in the nslet equation (\ref{GE:NS}) 
$u_{ki,0}+u_{{\bf k}j} u_{{\bf k}i,j}+p_{k,i}-\nu u_{ki,jj}= - \delta
(x_I-x_I') \delta _{ki}$ 
near to the origin of the nslet.
The delta function has order O$(L^{-3}T^{-1})$ in space $L$ and time $T$.
To obtain a non-trivial solution then the first term $u_{ki,0}$ has
the same order O$(L^{-3}T^{-1})$, which then means that the velocity
$u_{ki}$ has order O$(L^{-3})$.
This means that as we approach the origin of the nslet $L \rightarrow
0$ the term $u_{{\bf k}j} u_{{\bf k}i,j}$ is of order O$(L^{-7})$
and dominates the term $\nu u_{ki,jj}$ which is of order O$(\nu L^{-5})$
and so to leading order the nslet equation approaches the eulerlet
equation
\begin{equation}
  u_{ki,0}^E+u^E_{{\bf k}j} u^E_{{\bf k}i,j}+(p^B_{k,i}+p^E_{k,i})= - \delta
(x_I-x_I') \delta _{ki} ,
\label{eulerlet}
\end{equation}
where the pressure is decomposed into a term $p^B_k$ symmetric in $k$
and a term $p^E_k$ antisymmetric in $k$ from the following argument.
The term $\delta _{ki}$ suggests terms equating to it that are
antisymmetric in $k$ and $i$ so give zero integral contribution when
$i \ne k$ and a contribution when $i=k$.
However, the term $u^E_{{\bf k}j} u^E_{{\bf k}i,j}$ is symmetric in
$k$.
So this leads to the decomposition of the pressure as given,
and the eulerlet equation (\ref{eulerlet}) such that 
\begin{equation}
u^E_{ki,0}+p^E_{k,i} = - \delta \delta _{ki}  
\label{ER:eulerlet}
\end{equation}
where $p^E$ is the eulerlet pressure, and 
\begin{equation}
u^E_{{\bf k}j} u^E_{{\bf k}i,j}+p^{B}_{k,i}= 0
\label{ER:eulerlet_firstorder}
\end{equation}
where $p^B$ is the Bernoulli pressure \cite{Chadwick:2019}.

A solution to the eulerlet equation is now found by first
differentiating through (\ref{ER:eulerlet}) with respect to $i$ noting
that from the continuity equation then $u^E_{ki,i}= 0$ holds to give
\begin{equation}
  p^E_k =\frac{\delta(x_0-x_0')}{4\pi} \left[ \frac{1}{R} \right]
  _{,k} , 
\label{ER:stokeslet_pressure}
\end{equation}
since $\delta = \delta (x_0) \delta ({\bf x})$ and $\delta ({\bf x}) =
[-1/4\pi R]_{,ii}$.
The radius $R$ is the radial distance measure from a co-ordinate
system that is a moving reference frame.
In particular, within the boundary integral formulation the reference
frame origin centred on the eulerlet must move with the fluid velocity
at that point to ensure that the velocity vector direction is uniquely
defined.
Putting the expression for the pressure back into (\ref{ER:eulerlet})
and integrating gives
\begin{equation}
  u^E_{ki} = - H(x_0-x_0') \delta ({\bf x^L}) \delta _{ki} -
  \frac{1}{4\pi} 
  H(x_0-x_0') \left[ \frac{1}{R} \right] _{,ki}
\label{ER:eulerlet-velocity}
\end{equation}
where ${\bf x^L}$ is the distance from the moving 
co-ordinate reference frame, 
and $H(x_0-x_0')$ is the Heaviside function $H(x_0-x_0') =1 $ for
$x_0 > x_0'$ and zero otherwise.

%To see a description of the eulerlet potential $\Phi ^E = 1/R$ and the
%eulerlet velocity in (\ref{ER:eulerlet-velocity}),
%refer to the figures \ref{pot_S} and \ref{velfield_S}
%for the stokeslet potential and velocity,
%which are little different in form for the particular parameter values
%chosen there. 
For $R>0$, then $u^E_{ki} = \phi ^E _{k,i} = - \frac{1}{4\pi}
\left[ \frac{1}{R} \right] _{,ki}$, and
\begin{eqnarray}
  p^B_{k,i} &=& -\phi ^E _{{\bf k},j}\phi ^E _{{\bf k},ji} 
  = \left[ -(1/2) \phi ^E _{{\bf k},j}\phi ^E _{{\bf k},j} \right]_{,i} \nonumber \\
  p^B_k &=& -(1/2)  \phi ^E _{{\bf k},j}\phi ^E _{{\bf k},j}  = -(1/2)  u^E _{{\bf k}j}u^E _{{\bf k}j}  ,
\end{eqnarray}
which is the Bernoulli pressure.

\begin{figure}[h]
\begin{center}
\includegraphics[width=8.cm]{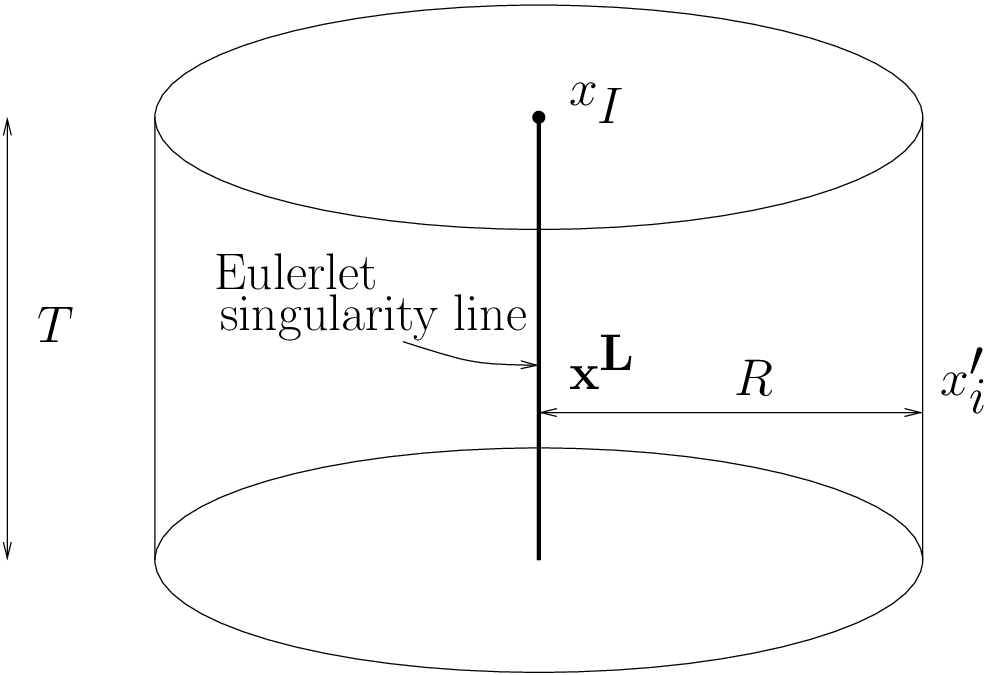}
\end{center}
\caption{The boundary $\partial \Sigma _{\delta}$ about the fluid
  point $x_I$ presented in a reference frame moving with the fluid at
  that point.} 
\label{movrefdelta}
\end{figure}

We can check the eulerlet contribution from (\ref{ER:eulerlet}) around
the boundary $\partial \Sigma_{\delta}$ consisting of the set of
points $x_i'$ centred on $x_I$ including: 
those such that $(x_j^L-{x_j}')(x_j^L-{x_j}')=R^2$ and $0 \le x_0-x_0' \le T$;
those such that $0 \le (x_j^L-{x_j}')(x_j^L-{x_j}') \le R^2$ and
$x_0-x_0'=T$, see figure \ref{movrefdelta}.

So, $\partial \Sigma_{\delta}$ and all the boundary surfaces
considered are piecewise continuous Lyapunov surfaces so have defined
normals up to the piecewise corners,
and so the Green's boundary integral and flux conditions hold.
The values of $T$ and $R$  are considered small,
the volume $V_R$ is the volume of the sphere of radius $R$,
and the surface $S_R$ is the surface of the sphere of radius $R$.
So, when a space-time infinitesimal element $d\sigma' $ is
over the same time, then $d\sigma' =dV'$, the volume infinitesimal
element,
and when a space-time infinitesimal element $d\sigma $
moves with the fluid, then $d\sigma' =ds' dt'$, where $ds'$ is an
infinitesimal surface and $dt'$ is an infinitesimal change in time.

So the eulerlet contribution from (\ref{ER:eulerlet}) around the
boundary $\partial \Sigma_{\delta}$ is 
\begin{eqnarray}
  \int _{\partial \Sigma _{\delta}}
    (u^E_{ki}n_0 +p^E_k n_i) d\sigma' &=&
  \int _{V_R} \left( -\delta (x_i^L-{x_i}') \delta _{ki} - \frac{1}{4\pi} \left[ \frac{1}{R}
    \right] _{,ki} \right) dV' \nonumber \\ 
  &&    + \int _T \int _{S_R} \frac{\delta (x_0-x_0')}{4\pi } \left[
    \frac{1}{R} \right] _{,k}  n_i ds' dt' \nonumber \\
    &=& -\delta _{ki} -\frac{1}{4\pi} \int _{V_R} \left[ \frac{1}{R} \right]_{,ki} dV' \nonumber \\
&&    - \frac{1}{4\pi } \int _{S_R} \frac{(x_k^L-{x_k}')}{R^3} \frac{(x_i^L-{x_i}')}{R} ds' \nonumber \\
    &=& -\delta _{ki} + \frac{1}{3} \delta _{ki} - \frac{1}{3} \delta _{ki} \nonumber \\
    &=& -\delta _{ki}.
\label{ER:check}
\end{eqnarray}

Additionally, it is noted that the viscous contribution in
(\ref{GE:NS}) is zero since
\begin{eqnarray}
&&  \int _{\partial \Sigma _{\delta}} \nu u^E_{ki,j}n_j d\sigma  \nonumber \\
&=&  -\frac{\nu}{4\pi} \int _0^T \int _{S_R} \left[ \frac{1}{R}
    \right]_{,kij}n_j ds dt \nonumber \\
  &=& -\frac{\nu T}{4\pi} \int _{S_R} \left[ -\frac{(x_k^L-{x_k}')}{R^3}
    \right]_{,ij}n_j ds \nonumber \\
&=& -\frac{\nu T}{4\pi} \int _{S_R} \left[ \frac{R^3\delta _{ki}-3R (x_i^L-{x_i}')(x_k^L-{x_k}')}{R^6}
    \right]_{,j}n_j ds \nonumber \\
&=& -\frac{\nu T}{4\pi} \int _{S_R} \left[ \frac{-3R^2\delta _{ki}+9(x_i^L-{x_i}')(x_k^L-{x_k}')}{R^6}
    \right] ds \nonumber \\
&=& -\frac{\nu T}{4\pi} \int _{S_R} \left[ \frac{-3R^2\delta _{ki}+9\delta
      _{ki} (x_1^L-{x_1}')(x_1^L-{x_1}')}{R^6}
    \right] ds \nonumber \\
&=& -\frac{\nu T}{4\pi} \int _{S_R} \left[ \frac{-3R^2\delta _{ki}+3(x_j^L-{x_j}')(x_j^L-{x_j}')\delta _{ki}}{R^6}
    \right] ds \nonumber \\
  &=& 0,
  \label{ER:viscous}
\end{eqnarray}
where we have used $x^L_{k,i}=\delta _{ki}$, $(x_j^L-{x_j}')(x_j^L-{x_j}')=R^2$,
and we have made some use of the symmetry properties of the integral.
So, the eulerlet also satisfies the equation for the nslet
(\ref{GE:NS}).
However, the eulerlet is not the nslet even though it satisfies the
nslet equation because it does not decay in the far-field.

Similarly, we can show equation (\ref{ER:eulerlet_firstorder}) holds
since
\begin{eqnarray}
&&  \int _{\partial \Sigma _{\delta}} u^E_{{\bf k}j} u^E_{{\bf
      k}i}n_j+p^{B}_kn_id\sigma \nonumber \\
&=&  \int _{\partial \Sigma _{\delta}} \phi^E_{{\bf k},j} \phi^E_{{\bf  k},i}n_j
  -(1/2)\phi^E_{{\bf k},j} \phi^E_{{\bf  k},j}n_i d\sigma \nonumber \\
  &=&0
  \label{quadraticsymm}
\end{eqnarray}
from the symmetry of the boundary ${\partial \Sigma _{\delta}}$.

\subsection{Far from the origin: the stokeslet}
Far from the origin, the velocity decays to zero, and so the quadratic
term disappears giving the time-dependent stokeslet defined by 
\begin{equation}
u^S_{ki,0}+p^S_{k,i} -\nu u^S _{ki,jj}= - \delta \delta _{ki}  \; , \;
\; u^S_{ki,i} = 0 .
\label{stokeslet}
\end{equation}
We note that this is also the time dependent oseenlet in a reference
frame moving with velocity $U$,
more usually attributed to far-field representations.
The steady stokeslet solution cannot be used here, as it does not
satisfy the unit force impulse but instead the unit force, as the term
$\delta $ is the Dirac delta function in four-dimensional space time
and not in three-dimensional space.
This fundamental solution for the Stokes equation appears to have been
given first by Oseen \cite{Oseen:1927} and later named the stokeslet.
The stokeslet solution is given in the appendix \ref{stoksletderiv}
from Fourier Transforms following the approach of Chan and Chwang
\cite{Chwang:2000} as 
\begin{equation}
  u^S_{ki} = \frac{1}{4\pi \nu} H(x_0-x_0')
  \left( \frac{{\rm erf} \eta }{R} \right) _{,0} \delta _{ki} 
  - \frac{1}{4\pi} H(x_0-x_0')
  \left( \frac{{\rm erf} \eta }{R} \right) _{,ki}
\label{stokeslet-velocity}
\end{equation}
where $\eta = \frac{R}{\sqrt{4\nu (x_0-x_0')}}$, 
and 
\begin{equation}
  p^S_k =\frac{\delta(x_0-x_0')}{4\pi} \left[ \frac{1}{R} \right] _{,k} 
\label{stokeslet-pressure}
\end{equation}
for the pressure.
%The potential $\Phi ^S = {\rm erf} \eta /R$ which generates this flow is
%represented in the following figure \ref{pot_S}.
%\begin{figure}[h]
%\begin{center}
%\includegraphics[width=9.cm]{potential_S.eps}
%\end{center}
%\caption{Representation of the potential $\Phi ^S$ with $4\pi \nu =1$}
%\label{pot_S}
%\end{figure}
%Considering the force impulse  with $k=1$ in (\ref{stokeslet-velocity})),
%this potential then generates a velocity field represented in the following
%figure \ref{velfield_S}.
%\begin{figure}[h]
%\begin{center}
%\includegraphics[width=9.cm]{velfield_S.eps}
%\end{center}
%\caption{The velocity field described by a vector plot}  
%\label{velfield_S}
%\end{figure}
%Here, we have taken $4\pi \nu =1$, $x_0=1$ and $x_3=x_3'$.
%This impulse force is in the $-x_1$ direction and applied at the nslet
%origin $x_1=x_1'$ generating a velocity field around it, see (\ref{GE:NS}).
%So as expected from figure \ref{velfield_S}, it is seen that
%there is a strong fluid velocity in the $x_1$ direction close
%to the nslet origin since from Newton's third law the fluid acts to
%apply an equal and opposite force at the origin point in the $x_1$
%direction. 
%Also interesting to note is that the vector directions trace out
%streamlines similar to that of a dipole.
%Due to the symmetric nature of the potential described in figure
%\ref{pot_S},
%the same descriptions are given for the other stokeslet approximations
%of the nslets, the only difference being that they are directed along
%the other axes. 

\section{Formulation of the Boundary Integral}
We now follow Green's Boundary Integral formulation in
\cite{Chadwick:2019a}, 
as used by Oseen \cite{Oseen:1927}.
Oseen starts with the integral given in the theory of hydrodynamic
potentials (described, for example, in \cite{Ladyzhenskaya:1969}) 
\begin{equation}
  \int _{\Sigma} \left[ f_i u_{ki} - f_{ki}u_i \right] d\Sigma ' = 0
\label{BI:Oseen}
\end{equation}
integrated over the primed Cartesian co-ordinate $x'_I$ such that in
the integrand the terms are functions of the following variables:
$f_i(x'_I)$ given by equation (\ref{GE:NavStokes}), $u_i(x'_I)$ which
is the fluid velocity, 
$u_{ki}(x_I-x_I')$ which is the $ith$ component of velocity for the
$kth$ nslet Green's function centred at $x_I=x_I'$,
and $f_{ki}(x_i-x_I')$ which is given by equation (\ref{GE:NS}).

\begin{figure}[h]
\begin{center}
\includegraphics[width=12.cm]{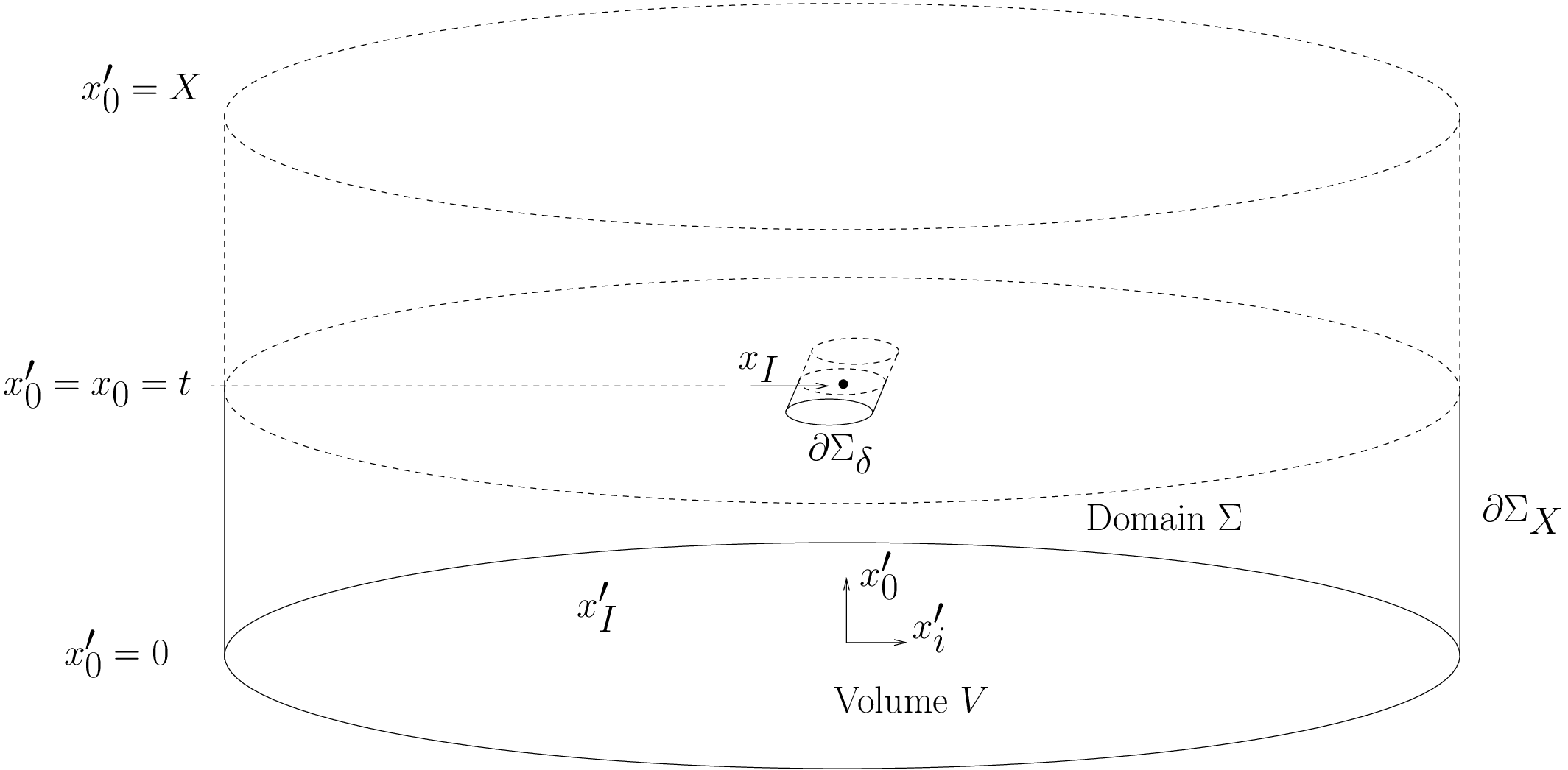}
\end{center}
\caption{The space-time domain $\Sigma$}   
\label{BI:domain}
\end{figure}

Since the nslet is undefined at its origin $x_I=x'_I$,
then this point is omitted from $\Sigma $.
Define $\Sigma $ as the domain between the boundaries $\partial
\Sigma_{\delta}$ and $\partial \Sigma_{X}$.
Here, the boundary $\partial \Sigma_{\delta}$ is defined in section
\ref{section4a}.
The boundary $\partial \Sigma_X$ is defined such that
it consists of the set of points $x'_I$ including:
those such that $(x'_j-x_j)(x'_j-x_j)=X^2$ and $-X \le x'_0 < 0$ or
$0 < x'_0 \le X$,  
those such that $0 \le (x'_j-x_j)(x'_j-x_j) \le X^2$ and $x'_0=0,X$,
for large $X$.
The domain $\Sigma $ is represented by the figure \ref{BI:domain}. 

%In the limit, no momentum flux escapes from this domain either from
%the fluid or the nslet because the far-field boundary of the domain
%$\partial \Sigma _X$ is such that the velocities all tend to zero,
%and the near-field domain $\partial \Sigma _{\delta}$ has been defined such
%that the fluid velocity and nslet direction are aligned to it.
%So the contributions from the quadratic terms $u_ju_{i,j}$ and
%$u_{{\bf k}j}u_{{\bf k}i,j}$ that arise in the Euler co-ordinate
%description of this flux must be zero. 
%Therefore it is expected that in the limit
%\begin{equation}
%  \lim _{X-\rightarrow \infty, \delta \rightarrow 0}
%  \left\{
%  \int _{\Sigma _X-\Sigma _{\delta}}
%  [ f_i^* u_{ki}-f_{ki}^*u_i] d\Sigma 
%  \right\}
%  =0
%\label{LagrangeBIform}
%\end{equation}
%holds, where $f_i^*=u_{i,0}+p_{,i}-\nu u_{i,jj} $ and $f_{ki}^*
%=u_{ki,0}+p_{k,i} -\nu u_{ki,jj}$ don't contain any non-linear terms
%$u_ju_{i,j}$ or $u_{{\bf k}j}u_{{\bf k}i,j} $.
%However, to be confident of this expectation,
%the non-linear contributions in the Euler co-ordinate description 
%are evaluated explictly in this section,
%so $f_i$ and $f_{ki}$ in (\ref{BI:Oseen}) are evaluated rather than
%$f_i^*$ and $f_{ki}^*$ in (\ref{LagrangeBIform}),
%and it works out that the non-linear terms give no contribution to
%the integral evaluation as expected.

On the domain boundary $\partial \Sigma _{\delta}$ that moves with the
fluid,
aligning the nslet with the fluid direction means that the quadratic
term in the nslet equation (which becomes the eulerlet equation in the
limit) drops out to leading order due to the symmetry of the eulerlet
(\ref{quadraticsymm}) in the integration calculated using a reference
frame moving with the fluid.
Similarly, the quadratic term in the Navier-Stokes equation drops out
to leading order because the material derivative approaches the
partial derivative since the boundary remains stationary in the moving
reference frame in the limit approaching the fluid point.  
So, the Lagrangian and Eulerian descriptions become the same, and so
the quadratic term vanishes which originates from the non-linear term
in the Navier-Stokes equation that appears as a consequence of the
material time derivative in the Eulerian flow description (and the
non-linearity doesn't appear for the material time derivative in the
Lagrangian flow description). 
So, the non-linearity in the boundary integral formulation is removed
allowing standard boundary integral theory to be used to formulate a
boundary integral representation of the velocity.

However, to be confident of this,
leading and lower order calculations are evaluated in detail in the
asymptotic limits next. 
Following \cite{Chadwick:2019a}, using the fact that
$[g(x_I-x'_I)]_{,I'} = -[g(x_I-x'_I)]_{,I}$, for a function $g$,
the integrand is rewritten as

\begin{eqnarray}
  f_iu_{ki}-f_{ki}u_i &=&
  \left\{ u_{i,0'}+u_ju_{i,j'}+p_{,i'}-\nu u_{i,j'j'} \right\} u_{ki}
  \nonumber \\
&&-\left\{ u_{ki,0}+u_{{\bf k}j}u_{{\bf k}i,j}+p_{k,i}-\nu u_{ki,jj}
  \right\} u_i
  \nonumber \\
&=& \left\{ u_{i,0'}+u_ju_{i,j'}+p_{,i'}-\nu u_{i,j'j'} \right\} u_{ki}
  \nonumber \\
&&-\left\{ -u_{ki,0'}-u_{{\bf k}j}u_{{\bf k}i,j'}-p_{k,i}-\nu u_{ki,j'j'}
  \right\} u_i
  \nonumber \\
  &=& (u_iu_{ki})_{,0'}+(pu_{ki})_{,i'}+(p_ku_i)_{,i'} \nonumber \\
  &&+(\nu u_{ki,j'}u_i-\nu u_{i,j'}u_{ki})_{,j'} \nonumber \\
&&+u_ju_{i,j'}u_{ki}+(u_{{\bf k}j}u_{{\bf k}i}u_i)_{,j'}-u_{i,j'}u_{{\bf
    k}i}u_{{\bf k}j} \nonumber \\
  &=& \left[ u_iu_{ki}\delta _{0J}+pu_{kj}+p_ku_j \right. \nonumber \\
&&\left. +\nu u_{ki,j'}u_i-\nu u_{i,j'}u_{ki} +u_{{\bf k}j}u_{{\bf
        k}i}u_i \right]_{,J'}
  \nonumber \\
&&+u_ju_{i,j'}u_{ki}-u_{i,j'}u_{{\bf k}i}u_{{\bf k}j} \nonumber \\
&=& \left[ u_i f_{kiJ}-u_{ki}\tau _{ij} \right] _{,J'} +Q_k,
\label{BI:Integrand}
\end{eqnarray}
where
$\tau _{ij} = -p\delta _{ij}+\nu u_{i,j'}$ is the stress tensor,
and 
$Q_k$ is the quadratic term contribution (in the sense that it
originates from the quadratic term present in the differential equations
(\ref{GE:NavStokes}) and (\ref{GE:NS})) defined by

\begin{eqnarray}
Q_k &=& u_ju_{i,j'}u_{ki}-u_{i,j'}u_{{\bf k}i}u_{{\bf k}j} \nonumber \\
&=& u_ju_{i,j'}u_{ki}-[u_iu_{{\bf k}j}u_{{\bf k}i}]_{,j'}+u_iu_{{\bf k}j}u_{{\bf k}i,j'} \nonumber \\
&=& u_ju_{i,j'}u_{ki}+u_iu_{{\bf k}j}(u_{{\bf k}i,j'}-u_{{\bf k}j,i'}) \nonumber \\
&&-[u_iu_{{\bf k}j}u_{{\bf k}i}]_{,j'}+[(1/2)u_iu_{{\bf k}j}u_{{\bf k}j}]_{,i'} \nonumber \\
&=& u_ju_{i,j'}u_{ki}+u_iu_{{\bf k}j}(u_{{\bf k}i,j'}-u_{{\bf k}j,i'}) \nonumber \\
&&-[u_iu_{{\bf k}j}u_{{\bf k}i}-(1/2)u_ju_{{\bf k}i}u_{{\bf k}i}]_{,j'},
\label{BI:Q}
\end{eqnarray}

\subsection{The non-linear divergence}
\begin{figure}[h]
\begin{center}
\includegraphics[width=12.cm]{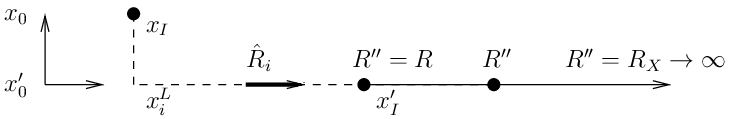}
\end{center}
\caption{The integral path of the non-linear potential in the moving
  reference frame about the fluid point $x_I$}
\label{BI:spoke}
\end{figure}

The quadratic term $Q_k$ is then given in terms of the divergence of a
vector field $q_{kj}$ directed on a radial spoke where
\begin{equation}
Q_k = q_{kj,j'}
\end{equation}
and 
\begin{equation}
  q_{kj}=W_k \hat{R}_j+q^*_{kj}
\label{qkj}
\end{equation}
such that $\hat{R}_j=\frac{{x_j}'-x_j^L}{R}$ where $R=\lvert
{x_j}'-x_j^L \rvert $ and
$q^*_{kj}=(1/2)u_ju_{{\bf k}i}u_{{\bf k}i}-u_iu_{{\bf k}j}u_{{\bf k}i}$.
Therefore the tensor $W_k \hat{R}_j$ has value $W_k$ directed along
the radial spoke $\hat{R}_j$, see figure \ref{BI:spoke},
and so
\begin{eqnarray}
  [W_k\hat{R}_j]_{,j'} 
  &=& u_ju_{i,j'}u_{ki}+u_iu_{{\bf k}j}(u_{{\bf k}i,j'}-u_{{\bf k}j,i'})
  \nonumber \\
  &=& g_k(x_I,x'_I).
\label{gk}
\end{eqnarray}
This is the divergence of a four dimensional vector in space-time,
and so expressing in the co-ordinate system $(x_0,R,\theta ,\alpha)$
where $(R,\theta ,\alpha)$ are spherical co-ordinates in space,
gives
\begin{equation}
  [W_k\hat{R}_j]_{,j'} = \frac{1}{R^2} \frac{\partial }{\partial R}
  \left\{ R^2W_k \right\} = g_k.
\end{equation}
So, integrating gives
\begin{equation}
  W_k({\bf x},{\bf x'}) = -\frac{1}{R^2} \int _{R_X}^{R}
  R''^2 g_k({\bf x},{\bf x''}) dR'',
  \label{Wk}
\end{equation}
where $R_X$ is the radial distance to the outer boundary $\partial
\Sigma _X$ of the domain $\Sigma $,  
the variable of integration is $R''=\lvert {x_i}''-x_i\rvert $,
such that ${x_i}''-x_i=R''\hat{R}_i$, $x_0=x_0'$,
and when $R''=R$ then ${x_i}''={x_i}'$, $x_0=x_0'$, 
see figure \ref{BI:spoke}.  

So, by putting the results from (\ref{BI:Integrand}) and
(\ref{BI:Q}) into (\ref{BI:Oseen}),
we can apply the divergence theorem to get a boundary integral
formulation 
\begin{eqnarray}
  \int _{\Sigma}
  \left[ f_i u_{ki} - f_{ki}u_i \right] d\Sigma ' &=&
  \int _{\Sigma}
  \left[ u_i f_{kiJ}-u_{ki}\tau _{ij} +q_{kj} \right] _{,J'} d\Sigma '
\nonumber \\
&=&   \int _{\partial \Sigma}
  \left[ u_i f_{kiJ}-u_{ki}\tau _{ij} +q_{kj} \right] n_J^{\Sigma} d\sigma ',
  \nonumber \\
  \label{BI:formulation}
\end{eqnarray}
where $n_J^{\Sigma}$ is the outward pointing normal to $\Sigma$,
and $\partial \Sigma $ is the collection of boundaries enclosing both
the exterior and interior domains that constitute $\Sigma $.
The first term in the integrand yields an integral associated with the
double layer potential,
the second term in the integrand with the single layer potential,
see for example the theory of hydrodynamic potentials given in
Ladyzhenskaya \cite{Ladyzhenskaya:1969}. 
The third term in the integrand is given by the non-linear potential
that depends on integration of quantities within the whole domain.
The ansatz proof assumes that the physical variables are all smooth,
and so the divergence theorem applied to the non-linear term to get
the boundary representation is possible since the term is assumed to
be smooth.
This representation also has similarities with the dual-reciprocity
method used for boundary integral representations for non-linear
partial differential equations,
except that the dual-reciprocity method approximates the non-linear
domain term by a series expansion of known functions,
such as radial basis functions,
whose divergence is known and so can then be represented over the
boundary.  

The whole concept of the boundary integral framework requires
linearity,
and the non-linear third term in (\ref{BI:formulation}) appears to
preclude this.
However, in the subsequent analysis we shall see that a representation
is found such that the non-linear term (as well as the double-layer
potential) does not give any singularities particularly for the
evaluation of the hypersurface boundary enclosing the fluid point in
the limit as the point is approached.
This is because part of the singularity (in the limit) from the
non-linear term in the velocity of the fundamental solution cancels
with the singularity (in the limit) of the eulerlet pressure.
The other part of the singularity (in the limit) gives a vorticity
term factor which is zero in the limit as the fundamental solution
approaches the eulerlet velocity.
This analysis means that the non-linear term does not
give a singularity in the limit, and instead is able to be evaluated
in the same way as for linear differential equations.  

\section{Evaluation of the Boundary Integral}\label{section6}

\subsection{Boundary integral evaluation in the vicinity of the point}
A similar approach to Chadwick \cite{Chadwick:2019a} is given, but
here in greater detail. 
We calculate part of the integral in
(\ref{BI:formulation}) which is around the integral path
  ${\partial \Sigma _{\delta}}$ in the vicinity of the point,
  since here the nslet approximates to the eulerlet.
From (\ref{BI:formulation}), the contribution to this integral from
the integral path ${\partial \Sigma _{\delta}}$ is
\begin{eqnarray}
I_{\delta} &=&  \int _{\partial \Sigma _{\delta}}
  \left[ u_i f_{kiJ}-u_{ki}\tau _{ij}+q_{kj} \right] n^{\Sigma}_J
  d\sigma ' \nonumber \\
  &=&
  -\int _{\partial \Sigma _{\delta}}
    \left[ u_i f_{kiJ}-u_{ki}\tau _{ij}+q_{kj} \right] n_J
  d\sigma ' \nonumber \\
  &=&
  -\int _{\partial \Sigma _{\delta}}
  \left[ u_i (u_{ki}\delta _{0J}+ u_{{\bf k}j}u_{{\bf k}i}+p_k
    \delta _{ij}-\nu u_{ki,j}) \right. \nonumber \\
    && \left. -u_{ki}\tau _{ij}+ (1/2)u_ju_{{\bf k}i}u_{{\bf
        k}i}-u_iu_{{\bf k}j}u_{{\bf k}i} + W_k \hat{R}_j \right] n_J  
  d\sigma ' \nonumber \\
  &=&
  -\int _{\partial \Sigma _{\delta}}
  \left[ u_i u_{ki}\delta _{0J} +u_jp_k -\nu u_i u_{ki,j} \right. \nonumber \\
    && \left. -u_{ki}\tau _{ij}+ (1/2)u_ju_{{\bf k}i}u_{{\bf
        k}i} + W_k \hat{R}_j \right] n_J  
  d\sigma '. \nonumber \\
\end{eqnarray}
Only the first term in the integrand involves the normal component
$n_0$,
all others have the factor $n_j$.
For $n_j$, then $R>0$ over the boundary and so
\begin{eqnarray}
  p_k &=& p_k^B+p_k^E \nonumber \\
&=&  -(1/2) u^E_{{\bf k}j}u^E_{{\bf k}j}+\frac{\delta (x_0)}{4\pi}
  \left[ \frac{1}{R} \right] _{,k}
\end{eqnarray}
where $u_{ki}$ tends towards $u^E_{ki}$.
So, the integral evaluation becomes
\begin{eqnarray}
  I_{\delta} &=& 
  -\int _{\partial \Sigma _{\delta}}
  \left[ u_i u^E_{ki}\delta _{0J} +u_jp^E_k -\nu u_i u^E_{ki,j} \right. \nonumber \\
    && \left. -u^E_{ki}\tau _{ij}+ W_k \hat{R}_j \right] n_J  
  d\sigma ', \nonumber \\
\label{ER:integral}
\end{eqnarray}

Evaluate the integral of each term in the integrand of 
(\ref{ER:integral}) in turn.
We see that from (\ref{ER:check}) and (\ref{primetodomain}),
the first two terms $u_i u^E_{ki}\delta _{0J} + u_jp^E_k$ give an
integral contribution $-u_k$.
Similarly, from (\ref{ER:viscous}),
the third term $-\nu u_i u^E_{ki,j}$ gives a zero leading order
integral contribution;
That is, the leading order in the Taylor series expansion of the
velocity
\begin{equation}
  u_i({\bf {x}'})=  u_i({\bf x^L})+u_{i,q}({\bf x^L})({x_q}'-x_q^L)
  +u_{i,ql}({\bf x^L})\frac{({x_q}'-x_q^L)({x_l}'-x_l^L)}{2!}+..
\end{equation}
The Taylor series expansion can be used because we are looking for a
solution for the class of variables in $C^{\infty}$ for which the
Taylor series can be employed.
The next order in the expansion gives a contribution
\begin{eqnarray}
  &&  \nu u_{i,q}\int _{\partial \Sigma _{\delta}}
u_{ki,j}x_q n_j d\sigma  \nonumber \\
  &=& -\nu \frac{u_{i,q}T}{4\pi} \int _S
  \left\{ -\frac{3\delta _{ki}}{R^4}+9\frac{(x_i^L-{x_i}')(x_k^L-{x_k}')}{R^6} \right\}
  (x_q^L-{x_q}') ds \nonumber \\
  &=&0
\end{eqnarray}
from symmetry.
The subsequent orders in the Taylor series expansion are all at most
of order O$(T) \rightarrow 0$. 

The next term in the integrand of (\ref{ER:integral}) is
$-u^E_{ki}\tau _{ij}$ and gives an integral contribution
\begin{eqnarray}
  &&
  \int _{\partial \Sigma _{\delta}} u^E _{ki} \tau _{ij} n_J d\sigma '
\nonumber \\
  &=&
  \tau _{ij} \int _{\partial \Sigma _{\delta}} u^E _{ki} n_J d\sigma '
  \nonumber \\
  &=& 
  -\tau _{ij} \int _{\partial \Sigma _{\delta}} u^E _{ki} n_J d\sigma 
  \nonumber \\
  &=& 
  -\tau _{ij} \int _0^T \int _S u^E _{ki} n_J ds dt
  \nonumber \\
  &=& 
  \frac{\tau _{ij}T}{4\pi} \int _S \left[ \frac{1}{R} \right]_{,ki}
  n_J ds \nonumber \\
  &=& 
  \frac{\tau _{ij}T}{4\pi} \int _S \left[ \frac{R^3 \delta
      _{ki}-3R(x_i^L-{x_i}')(x_k^L-{x_k}')}{R^6} \right] \frac{(x_j^L-{x_j}')}{R} ds \nonumber \\
  &=& 0
\end{eqnarray}
since subsequent terms in the Taylor series expansion are of order
O$(T) \rightarrow 0$.

Finally, the last term is
\begin{equation}
  -\int _{\partial \Sigma _{\delta}} W_k \hat{R}_j n_J d\sigma '
  = 
  -\int _T \int _S W_k \hat{R}_j n_J d\sigma '.
\label{lastterm}
\end{equation}
%and $R>0$ on the surface $S$, so $u^E_{{\bf k}i,j}=u^E_{{\bf k}j,i}$
%and
%$q^E_k= \frac{1}{R^2} \int _{R_X}^R R^{*2}g_k^EdR^*$.
The integral for $W_k$ in (\ref{Wk}) is split into three parts as follows
\begin{equation}
  W_k= -\frac{1}{R^2}
  \left\{ \int _{R_X}^{R_L} + \int _{R_L}^{R_{\epsilon}} + \int _{R_{\epsilon}}^R \right\} 
  R''^2g_kdR'' = I_1+I_2+I_3,
  \label{threeintegrals}
\end{equation}
where $R_L$ is a sufficiently large constant such that for $R''>R_L$
then the nslet approximates to the stokeslet $u_{ki}\sim u_{ki}^S $,
and where $R_{\epsilon}$ is a sufficiently small constant such that for $R''<R_{\epsilon}$
then the nslet approximates to the eulerlet $u_{ki}\sim u_{ki}^E$.
Each of the three integrals are calculated in turn.

First, it is seen that for $R_X$ sufficiently large then
$u_{ki} \sim u_{ki}^S \sim \frac{1}{{R''}^3}$.
This means that from (\ref{gk})
$g_k \sim \frac{1}{{R''}^3} $ and so, from (\ref{threeintegrals}), $I_1
\sim (1/R^2) \ln R_X$. 

Second, it is seen that between the two constant limits $R_L$ and
$R_{\epsilon}$ then $u_{ki} \sim constant $ (bounded).
This means that from  (\ref{threeintegrals}), $I_2 \sim 1/R^2 $. 

Third, it is seen that for $R$ sufficiently small then
$u_{ki} \sim u_{ki}^E \sim \frac{1}{{R''}^3}$ over the spherical
surface radius $R>0$.
(There is only a contribution from the spherical surface because there
is no contribution when $J=0$, since from (\ref{lastterm}) $\hat{R}_j
n_J =\hat{R}_j n_j $.) 
For $R>0$ the eulerlet vorticity is zero,
giving $u^E_{{\bf k}i,j}=u^E_{{\bf k}j,i}$,  
which can be easily verified from (\ref{ER:eulerlet-velocity}).
This means that the quadratic singular term of the fundamental
solution in $g_k$ given by (\ref{gk}) has a factor to leading order
$u^E_{{\bf k}i,j}-u^E_{{\bf k}j,i}=0$,
and this is the reason that despite the non-linear quadratic
singularity in the expression the boundary integral can still be
evaluated.
From appendix \ref{nsletviscorder},
it is seen that the next lower order viscous term is exponentially
small with the factor $e^{-\eta ''^2} \rightarrow 0 $ for $\eta ''^2
=R''^2/(4\nu x_0) \rightarrow \infty$, which holds since $T=o(R^2/\nu)$. 
This means that from (\ref{gk}) that
$g_k \sim \frac{1}{{R''}^3} $ and so, from (\ref{threeintegrals}), $I_3
\sim (1/R^2) \ln R$.

%and in the Green's integral representation the velocity is expected to
%have the same decay as the Green's function,
%and so $u_i \sim u_i^S \sim \frac{1}{R^{*3}}$.

So, from (\ref{threeintegrals}), $I_1 \sim (1/R^2) \ln R_X$, $I_2 \sim 1/R^2$
and $I_3 \sim (1/R^2) \ln R$,
and so the integral evaluation of (\ref{lastterm}) where boundary
$\partial \Sigma _{\delta}$ has time interval $T$ and spherical radius
$R$ small is
\begin{equation}
  \int _{\partial \Sigma _{\delta}} W_k \hat{R}_j n_J d\sigma '
  \sim T Maximum( \lvert \ln R_X \rvert, \lvert \ln R \rvert).
\end{equation}
Choosing $T$ to be sufficiently small, of smaller order than
$Maximum(\lvert \ln R_X \rvert, \lvert \ln R \rvert)$ ensures that the
integral contribution disappears in the limit.
For example, we can choose $R_X\sim 1/R$, $T\sim R^3$, and as $R
\rightarrow 0$, then the integral contribution goes as $R^3 \ln R
\rightarrow 0$.

So, bringing all the terms together gives the integral contribution
\begin{equation}
  I_{\delta} = -u_k.
\label{ER:I-delta}
\end{equation}

\subsection{Boundary integral evaluation far from the point} 

In the far-field around $\partial \Sigma _X$, the time interval $T$
and radius $R$ are large,
the nslet tends towards the stokeslet (which is the oseenlet in a
moving reference frame), and this contribution is shown to be zero
next.
In the far-field, we have the contribution
$u_{ki} \sim u_{ki}^S \sim \frac{1}{{R''}^3}$,
and $u_i \sim u_i^S \sim \frac{1}{{R''}^3}$.
So $g_k \sim \frac{1}{{R''}^{10}}$, and therefore
$W_k \sim \frac{1}{R^2} \frac{1}{R^7}= \frac{1}{R^9}$.
So, the far-field integral is 

\begin{eqnarray}
  \lvert I_X \rvert &=& \lvert
  \int _{\partial \Sigma _X}
  \left[ u_i f_{kiJ}-u_{ki}\tau _{ij} +q_{kj} \right] n_J^{\Sigma}
  d\sigma '
  \rvert
  \nonumber \\
  &\sim & \lvert \int _{\partial \Sigma _X}
  \left[ u_i f^S_{kiJ}-u^S_{ki}\tau _{ij} +q^S_{kj} \right] n_J^{\Sigma}
  d\sigma ' \rvert \nonumber \\
\end{eqnarray}
However, $\int _{\partial \Sigma _X} u_i f^S_{kiJ}n^{\Sigma}_J d\sigma
' \sim R^3 (1/R^3)(1/R^2)$ for a leading order pressure term decay
$1/R^2$ and the time integral $T$ being chosen of order $R$ large.
Similarly, $\int _{\partial \Sigma _X} u_{ki}^S \tau_{ij}n^{\Sigma}_j
d\sigma ' \sim R^3 (1/R^3)(1/R^2) $.
The integral contribution $\int _{\partial \Sigma _X} q_{ki}^S n^{\Sigma}_j
d\sigma ' \sim R^3 (1/R^9)$.
(The contribution to $q_{kj}$ from the term $W_k \hat{R}_j$ in
(\ref{qkj}) is identically zero because $R=R_X$ and so $W_k$ is
identically zero from (\ref{Wk}), which means that $g_k=0$ from
(\ref{gk}).)    
So
\begin{equation}
  I_X \rightarrow 0.
\label{farintcalc}
\end{equation}

\section{Velocity representation}\label{VR}
We now give an nslet representation for the boundary integral
(\ref{BI:formulation}).
In the near-field around $\partial \Sigma _{\delta}$,
the nslet tends towards the eulerlet and gives a contribution $-u_k$,
see (\ref{ER:I-delta}).
In the far-field around $\partial \Sigma _X$,
the nslet tends towards the stokeslet and gives zero contribution.
Therefore,
putting these results into (\ref{BI:formulation})
then gives the velocity by the boundary integral representation
\begin{equation}
  u_k = \int _{\partial \Sigma }
  \left[ u_i f_{kiJ}-u_{ki}\tau _{ij} +q_{kj} \right] n_J^{\Sigma} d\sigma '
\label{intrepforvel}
\end{equation}
for some boundary $\partial \Sigma $.
The equation (\ref{intrepforvel}) demonstrates why a representation
for the velocity is valid here,
because the calculation in the vicinity of the point
(\ref{ER:I-delta}) gives the velocity, 
and far from the point (\ref{farintcalc}) gives zero.
%Further, by combining fluid motions either side of the boundary
%$\partial \Sigma $ (both external and internal flow together), it
%shall be shown that the non-linear term $q_{kj}$ is continuous across
%the boundary and so disappears as well.  

%Consider two types of boundary, the first generated by a moving body and
%the second stipulating values at an initial condition.

\subsection{Initial Value Problem boundary}\label{IVP}
Consider the variables stipulated at time $x_0 =0$.
The outward pointing normal is $n_I^{\Sigma} = - \hat{x}_0$,
where $\hat{x}_0$ is the unit vector in the time direction.
So, from (\ref{intrepforvel})
\begin{equation}
  u_k = \int _V u_i f_{ki0}n_0^{\Sigma} dV' = -\int _V u_i u_{ki} dV',
\label{MP:IVP}
\end{equation}
where the space-time boundary $\partial \Sigma$ is chosen as the
volume $V$.
It is noted that the representation (\ref{MP:IVP}) contains a single
layer potential only. 
We note that throughout the volume $V$ at time $x_0 =0$,
the velocity $u_i$ has been stipulated as an initial condition,
and so is known.
Substituting the null homogeneous boundary condition $u_i=0$ in the
integrand gives the velocity $u_k=0$ for all subsequent time.

\section{Existence and Smoothness}\label{ESU}

\begin{figure}[h]
\begin{center}
\includegraphics[width=12cm]{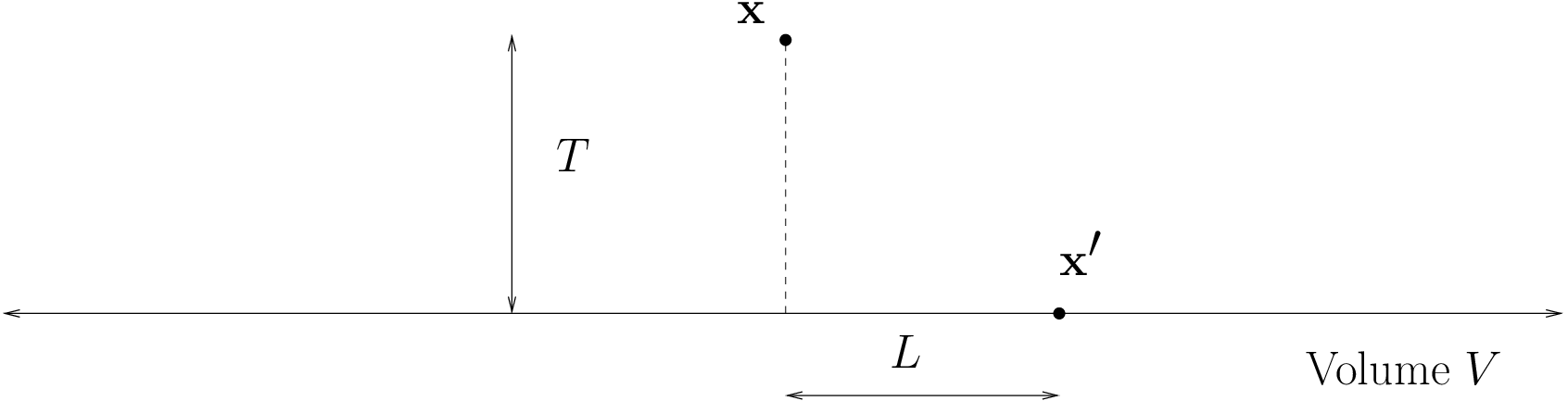}
\end{center}
\caption{The space-time representation of the volume integration for the velocity}   
\label{Volint}
\end{figure}

Putting in the function variables explicity into (\ref{MP:IVP}) gives
\begin{equation}
  u_k({\bf x},T) = -\int _V u_i({\bf x'},0) u_{ki}({\bf x - x'},T) dV',
\label{MP:IVPwithfunvar}
\end{equation}
with the integration shown diagrammatically by the space-time
representation in figure \ref{Volint}.

Consider $T$ sufficiently small.
Orders of terms in the nslet equation have been given in
\ref{section4a}, 
and following the same arguments then the nslet approximates to the
eulerlet to order $\epsilon $ such that
$u_{kj}\Delta V'=u_{kj}^E\Delta V'(1+$O$(\epsilon))$ where $\epsilon = \nu
L^{-5}/(L^{-3}T^{-1}) = \nu T/L^2$.
So as $L$ gets smaller such that ${\bf x'}$ becomes close to ${\bf x}$
then the error gets larger.
However, for $L$ sufficiently small it has also been given in
\ref{section4a} that the nslet approximates to the eulerlet to order
$\epsilon$ where $\epsilon = \nu L^{-5}/L^{-7} = \nu L^2$,
which is small for small $L$.
Making use of both bounds by choosing the first for $L > L_0$ and the
second for $L<L_0$ where $L_0 = T^{1/4}$ ensures that the nslet
approximates to the eulerlet to order $\epsilon = \nu \sqrt{T}$
throughout the whole range of the volume integration.

Therefore
\begin{eqnarray}
  u_k({\bf x^L},T) &=& -\int _V u_i({\bf x'},0)u_{ki}^E({\bf x-x'},T)
  dV' (1+{\rm O}(\nu \sqrt{T})) \nonumber \\
  &=& \left[ u_k({\bf x},0) + \frac{1}{4\pi} \int _V u_i({\bf x'},0)
    [1/R]_{,ki} dV' \right] (1+{\rm O}(\nu \sqrt{T})) 
  \nonumber \\
  &=& \left[ u_k({\bf x},0)- \frac{1}{4\pi} \int _{S_R} u_i({\bf x'},0)
    [1/R]_{,k} n_i ds' \right] (1+{\rm O}(\nu \sqrt{T}))   \nonumber \\
\end{eqnarray}
since $u_i[1/R]_{,ki}=-[u_i[1/R]_{,k}]_{,i'}$,
where ${\bf x^L}$ is the Lagrangian co-ordinate position moving with
the fluid, 
and the spherical surface $S_R$ has radius $R$ enclosing $V$.
However, on $S_R$ as $R \rightarrow \infty$ then $u_i \rightarrow 0$,
$S_R= $O$(R^2)$, $[1/R]_{,k}=$O$(R^{-2})$
and so the integral decays to zero in the limit leaving
\begin{equation}
  u_k({\bf x^L},T) = u_k({\bf x^L},0)(1+{\rm O}(\nu \sqrt{T})).
  \label{velocityorder}
\end{equation}
Similarly,
\begin{eqnarray}
  u_{k,j}({\bf x^L},T) 
  &=& u_{k,j}({\bf x},0)(1+{\rm O}(\nu \sqrt{T}))-\nonumber \\
  && \frac{1}{4\pi} \int _{S_R} u_i({\bf x'},0)
    [1/R]_{,kj} n_i ds' (1+{\rm O}(\nu \sqrt{T})) 
\end{eqnarray}
such that on $S_R$,
$[1/R]_{,kj}=$O$(R^{-3})$ and so the integral decays even faster to
zero giving
\begin{equation}
  u_{k,j}({\bf x^L},T) = u_{k,j}({\bf x^L},0)(1+{\rm O}(\nu \sqrt{T})).
  \label{velderivorder}
\end{equation}
Similarly, the same argument applies for all the other derviatives
meaning that the velocity is smooth.
From the appendix \ref{pressuresmooth}, this means that the pressure
is smooth also.

So the solution is smooth in the region
$\nu \sqrt{T} \le \epsilon$ for small $\epsilon$.
Choosing for example $\epsilon = 10^{-3}$, then this applies to a
timestep $T= \nu ^{-2}10^{-6}$.
Repeating for $N=10^6\nu ^2 t$ time steps, a general time $t$ is
reached.
So the solution exists and is smooth for a general time $t$.

This is a proof by ansatz, so it was originally assumed that the
velocity and pressure were smooth.
The formulation was then obtained and shown to represent a velocity
and pressure that are smooth for all time as originally assumed,
giving proof by ansatz.

\section{An example of chaotic motion} 

Consider the high Reynolds number limit of Euler flow.
The eulerlet has a singularity line extending out into the fluid.
This singularity line must be aligned to the fluid velocity so that
the velocity direction is uniquely defined.
The eulerlet singularity line is therefore directed along the
Lagrangian co-ordinate giving a dynamical system with the possibility
of chaotic mixing.
For example the blinking vortex \cite{DaitcheTel:2009},
is obtained in this formulation by considering a span-wise
distribution of fundamental solutions such that the force direction is
perpendicular to the flow field forward motion.
This is shown to reduce to a horseshoe vortex in \cite{Chadwick:2019}.
Close to a trailing vortex of this horseshoe,
the flow approximates to a two-dimensional steady point vortex.
It is possible to consider two such point vortices interacting with
each other such that their origins (centres) move with the fluid as a
direct consequence of the requirement that the singularity line moves
with the fluid. 
If we then suppose forces applied by an agitator such that the
strength of the vortices may vary as given in \cite{DaitcheTel:2009},
then the problem reduces to that of the blinking vortex in which chaos
is demonstrated. 
%Further, supposing we have high energy flow but not at the Euler
%limit, then shearing together with chaos gives turbulent motion in
%three-dimensional flow.

\section{Conclusion}

%A Green's integral representation for the velocity is given over a
%four-dimensional space-time hypersurface with a Green's function
%fundamental solution kernel, called an nslet,
%having strength given by a force impulse distribution over the
%hypersurface. 
%The force impulse for the problem,
%and so equivalently energy of the system,
%is assumed to be bounded.
A Green's integral representation for the velocty is given by a volume
integral at initial time with a Green's function fundamental solution
kernal, called an nslet, having strength given by (minus) the initial
velocity.
For a fluid point a small time step away,
then the nslet approximates to the eulerlet to within a known error
bound. 
%Consequently the velocity a small time step away from initial time is
%approximated by a velocity at initial time to within this known error.
%The same analysis applies to all the velocity derivatives as well, and
%the approach can be repeated via multiple time steps to reach any
%desired time.
Repeating for multiple time steps,
it is shown that the velocity is smooth and consequently also the
pressure at any later time. 
%close to the fluid point,
%then the nslet approximates to the eulerlet giving an
%integral representation by a known smooth function even though the
%eulerlet itself is singular in the domain.
%Hence this integral approach addresses and removes the possibility of
%finite-time singularities occuring
%(because point force impulses occur only on the space-time boundaries
%and not in the fluid)
%which has been a major difficulty in finding existence and smoothness
%proofs up until now.   
%From this, it is shown that if the velocity is smooth at an initial
%time then it is smooth for all subsequent time as well.
%If the velocity exists and is smooth then from the Helmholtz theorem
%so is the pressure.
The force impulse distribution is given by the initial velocity,
and if this is taken to be zero then the velocity at all later time is
zero.

The fundamental solution origin moves with the flow field,
allowing the possibility of a chaotic dynamical system and the example
of the blinking vortex is given.
Future work then is to model turbulence and test cascade down to the
expected Kolmogorov length scale.

\appendix
%\numberwithin{equation}{section}
%\begin{appendices}
  
\section{Stokeslet}\label{stoksletderiv}
We follow the approach of Chan and Chwang \cite{Chwang:2000} to give
the unsteady stokeslet.
The stokeslet equation is 
\begin{equation}
u^S_{ki,0}+p^S_{k,i} -\nu u^S _{ki,jj}= - \delta \delta _{ki}  \; , \;
\; u^S_{ki,i} = 0 .
\label{A:stokeslet}
\end{equation}
Taking the divergence of the first equation in (\ref{A:stokeslet})
gives $p^S_{k,ii}=-\delta _{,k}$.
Letting $p^S_k= \delta (x_0) p^*_{,k}$ gives $p^*_{,ii}=-\delta (x_i)$
and so $p^* = \frac{1}{4\pi R}$, which gives
\begin{equation}
  p^S_k =\frac{\delta(x_0)}{4\pi} \left[ \frac{1}{R} \right] _{,k} .
\label{A:stokeslet-pressure}
\end{equation}
Substituting this back into the first equation in (\ref{A:stokeslet})
gives
\begin{equation}
  u^S_{ki,0} -\nu u^S_{ki,jj}= -\delta \delta _{ki} - \left(
  \frac{1}{4\pi R} \right) _{,ki} \delta (x_0).
\end{equation}
We can split the velocity into two parts, and determine each
separately,
such that
\begin{equation}
  u^S_{ki} = u^A \delta _{ki} + u^B_{,ki}
\end{equation}
where
\begin{equation}
  u^A_{,0}-\nu u^A_{,jj} = -\delta
\label{A:uA}
\end{equation}
and
\begin{equation}
  u^B_{,0}-\nu u^B_{,jj} = -\frac{1}{4\pi R} \delta (x_0).
\label{A:uB}
\end{equation}

\subsection{Evaluation of $u^A$}
Apply Fourier Transforms to (\ref{A:uA}) to get
\begin{equation}
  \left( \frac{1}{\sqrt{2\pi}} \right)^4
    \int _{\Sigma _{\infty}}
    \left( u^A_{,0} -\nu u^A_{,jj} \right) e^{-\alpha _J x_J } d\Sigma
    =- \left( \frac{1}{\sqrt{2\pi}} \right)^4,
\end{equation}
where $\alpha _I$ are the Fourier Transform variables,
and $\int _{\Sigma _{\infty}}$ is the integration across all space-time.
So
\begin{equation}
  \left( \frac{1}{\sqrt{2\pi}} \right)^4
  \left( i\alpha _0 +\nu \alpha _R ^2 \right)
  \int _{\Sigma _{\infty}}   u^A
  e^{-\alpha _J x_J } d\Sigma
  = - \left( \frac{1}{\sqrt{2\pi}} \right)^4.
\end{equation}
Let ${\bar{\bar u}}^A =   \left( \frac{1}{\sqrt{2\pi}} \right)^4
\int _{\Sigma _{\infty}} u^A e^{-\alpha _J x_J } d\Sigma$,
and
let ${\bar u}^A =   \left( \frac{1}{\sqrt{2\pi}} \right)^3
\int _{V_{\infty}} u^A e^{-\alpha _j x_j } dV$,
where $dV$ is an element of volume $V_{\infty}$ which is an
integration across all of space.
This gives
\begin{eqnarray}
  {\bar {\bar u}}^A &=&
  - \left( \frac{1}{\sqrt{2\pi}} \right)^4
  \left( \frac{1}{i} \right) \left( \frac{1}{\alpha - i\nu \alpha
    _R^2} \right) \nonumber \\
  {\bar u}^A  &=&
    - \left( \frac{1}{\sqrt{2\pi}} \right)^3 H(x_0)
 e^{\nu \alpha _R^2 x_0}\nonumber \\
\end{eqnarray}
using the result for the inverse Fourier Transform \newline $F^{-1}\{ 1/(\alpha
-ia) \} = H(x_0) \sqrt{2\pi} i e^{-\alpha x}$ given in standard
transform tables,
for the real variable $x$, transform variable $\alpha $ and constant $a$.
Taking the inverse transforms in the spatial directions by using the
inverse Fourier Transform $F^{-1}\{ e^{-a\alpha ^2} \} = \left(
e^{-x^2/(4a)} \right) /\sqrt{2a}$ given in standard transform tables
for the Gaussian,
we then get
\begin{eqnarray}
  u^A &=& 
  - \left( \frac{1}{\sqrt{2\pi}} \right)^3
  H(x_0) e^{-R^2/(4\nu x_0)} 
  \left( \sqrt{\frac{1}{2x_0}} \right)^3 \nonumber \\
&=& \frac{1}{4\pi \nu} H(x_0) \left( \frac{{\rm erf} \eta}{R} \right)
    _{,0}
\end{eqnarray}
where $\eta = \frac{R}{\sqrt{4\nu x_0}}$ and 
${\rm erf} \eta = (2/\sqrt{\pi}) \int _0^{\eta} e^{-\eta '^2} d\eta
'$,
which means
\begin{eqnarray}
({\rm erf} \eta )_{,0} &=& \eta _{,0} ({\rm erf }\eta )_{,\eta } 
  \nonumber \\
  &=& -\frac{R}{4} \frac{1}{\sqrt{\nu}} \frac{1}{\sqrt{x_0}^3}
  \frac{2}{\sqrt{\pi}} e^{-\eta ^2} .
\end{eqnarray}

\subsection{Evaluation of $u^B$}
We then have
\begin{equation}
[u^B_{,0} - \nu u^B_{jj} ]_{,qq} = -\left( \frac{1}{4\pi R} \right)
_{,qq} \delta (x_0) = \delta .
\end{equation}
So $u^B_{,jj} = -u^A$ which gives $u^B_{,0} =-\nu u^A$ and so
\begin{equation}
  u^B = -\frac{1}{4\pi} H(x_0) \frac{{\rm erf} \eta }{R}.
\end{equation}

This then gives the stokeslet velocity to be
\begin{equation}
  u^S_{ki} = \frac{1}{4\pi \nu} H(x_0)
  \left( \frac{{\rm erf} \eta }{R} \right) _{,0} \delta _{ki} 
  - \frac{1}{4\pi} H(x_0)
  \left( \frac{{\rm erf} \eta }{R} \right) _{,ki}.
\end{equation}

\section{Pressure smooth}\label{pressuresmooth}
By applying the Helmholtz decomposition (fundamental theorem of vector
calculus) it is shown that if the velocity is smooth then the pressure
is smooth as well.
The Navier-Stokes equation is rewritten in terms of the head
function $h$ such that
\begin{equation}
  u_{i,0}+h_i+p_{,i}-\nu u_{i,jj} = 0 \; , \; \; h_i=u_ju_{i,j}.
\end{equation}
Applying the Helmholtz decomposition to the velocity and head
functions
\begin{eqnarray}
  u_i &=& \phi _{,i}+w_i \;, \; \; w_i = \epsilon _{ijk} A_{k,j}
\nonumber \\
    h_i &=& h_{,i}+h_i^w \;, \; \; h_i^w = \epsilon _{ijk} B_{k,j}
    \nonumber \\
\end{eqnarray}
for curl free scalar potentials $\phi$ and $h$,
and divergence free vectors $w_i$ and $h_i^w$,
represented by the vector potentials $A_k$ and $B_k$ respectively,
and $\epsilon _{ijk}$ is the Levi-Civita symbol.
(This symbol is such that $\epsilon _{ijk} = 1$ for permutations
$(i,j,k) = (1,2,3)$, $(2,3,1)$ and $(3,1,2)$,
$\epsilon _{ijk} = -1$ for permutations $(i,j,k) = (1,3,2)$, $(2,1,3)$ 
and $(3,2,1)$, and zero otherwise.)
The wake velocity divergence is $w_{i,i}= \epsilon
_{ijk}A_{k,ji}=0$, and so $u_{i,i} = \phi _{,ii} = 0$.
So, the curl free component of the Navier-Stokes equation is
$\phi _{,i0}+h_{,i} +p_{,i} = 0$ since $\nu \phi _{,ijj}=0$.
Integrating this then gives the pressure as
\begin{equation}
  p=-\phi _{,0}-h.
\end{equation}
From the Helmholtz theorem \cite{Koenigsberger:1906},
the potentials are given by
\begin{eqnarray}
  \phi &=& \frac{1}{4\pi} \int _S \frac{u_in_i}{\lvert R-R'\rvert } ds' \nonumber
  \\
  h&=& -\frac{1}{4\pi} \int _V \frac{u_{j,i'}u_{i,j'}}{\lvert R-R'\rvert } dV' +
  \frac{1}{4\pi} \int _S \frac{u_ju_{i,j'} n_i}{\lvert R-R'\rvert } ds', \nonumber \\
\end{eqnarray}
for some volume $V$ enclosed by a surface $S$ with outward pointing
normal $n_i$.
So, if the velocity $u_i$ is smooth then the potentials $\phi $ and $h$
are smooth, and the pressure $p$ is smooth.

\section{Next order viscous term for near-field nslet}\label{nsletviscorder}  
Consider the next order viscous term in the near-field expansion for
the nslet about the eulerlet, in particular when $R^2/(\nu T)
\rightarrow \infty$.

The leading order near-field expansion for the nslet is the eulerlet
given in \ref{section4a} to satisfy $u_{ki} \sim u_{ki}^E$ and $p_k
\sim p_k^B+p_k^E$ where $u_{ki,0}^E+p_k^E=-\delta \delta _{ki}$ and 
$u_{{\bf k}j}^Eu_{{\bf k}i,j}^E+p_{k,i}^B=0$.
So in the near-field such that $L, T \rightarrow 0$ then
\begin{equation}
u_{ki,0}+(p_k-p_k^B)_{,i}-\nu u_{ki,jj}=-\delta \delta _{ki}
-(p_{k,i}^B+u_{{\bf k}j}u_{{\bf k}i,j}),
\label{nsletexpvisc}
\end{equation}
such that $p_{k,i}^B+u_{{\bf k}j}u_{{\bf k}i,j}$ is small
on account of $p_{k,i}^B+u_{{\bf k}j}^Eu_{{\bf k}i,j}^E=0$.
Therefore, expand about this small term so
$u_{ki}=u_{ki}^S+u_{ki}^{II}+....$, and
$p_k-p_k^B=p_k^S+p_k^{II}+.....$
such that
\begin{eqnarray}
  u_{ki,0}^S+p_{k,i}^S-\nu u_{ki,jj}^S &=& -\delta \delta _{ki}
  \nonumber \\
  u_{ki,0}^{II}+p_{k,i}^{II}-\nu u_{ki,jj}^{II} &=& -p_{k,i}^B-u_{{\bf
      k}j}^Su_{{\bf k}i,j}^S. \nonumber \\
\end{eqnarray}
The first equation is the stokeslet, and
$-p_{k,i}^B=u_{{\bf k}j}^Eu_{{\bf k}i,j}^E$,  
whereby from the eulerlet (\ref{ER:eulerlet-velocity}) and
the stokeslet (\ref{stokeslet}), they approximate to order
$u_{kj}^S=u_{kj}^E(1+$O$(e^{-\eta ^2}))$
where $\eta \rightarrow \infty$
for $R/(\sqrt{\nu T}) \rightarrow \infty$.
Similarly, the second order term decays exponentially as well on
account of
$  u_{ki,0}^{II}+p_{k,i}^{II}-\nu u_{ki,jj}^{II} = {\rm O}(e^{-\eta
  ^2}/R^6)$ as exponential decay is faster than any polynomial.

%\end{appendices}

%\ethics{This article has no ethics issues.}
%\noindent {\bf Ethics} This article has no ethics issues.

%\dataccess{This article has the additional data which is gnu fortran
%  f95 source code for the new method, the shooting method and Kusukawa 
%  solution for uniform flow past a semi-infinite flat plate.}
%\noindent {\bf Data access} This article has the additional data which
%is gnu FORTRAN f95 source code for the new method, the shooting method
%and Kusukawa solution for uniform flow past a semi-infinite flat plate.

%\aucontribute{This is the sole work of EC.}
%\noindent {\bf Author contribution} This is the sole work of EC.  

%\competing{There are no competing interests.}
%\noindent {\bf Competing interests} There are no competing interests.

%\funding{The work on the verification (Section 11) has been funded by
%  EPSRC grant EP/V058754/1 entitled: Mathematical verification for a
%  novel Navier-Stokes representation.} 
%\noindent {\bf Funding} The work on the verification (Section 11) has
%been funded by EPSRC grant EP/V058754/1 entitled: Mathematical
%verification for a novel Navier-Stokes representation. 

%\ack{}

%\section*{Acknowledgement}
%The author wishes to acknowledge support from EPSRC through the funded
%grant above,
%and the research team Dr. Rabeea Darghoth, Dr. Bwebum Dang and
%Mr. Hamid Adamu.  

%\section*{References}
%\bibliographystyle{jfm}
\bibliographystyle{plain}
\bibliography{chadwick}

\end{document}